\documentclass[12pt,leqno]{article}
\usepackage[final]{graphics}
\usepackage{amsthm,amsfonts,epsfig}
 \usepackage{latexsym, amssymb}

\def\ds{\displaystyle}

\def\ni{\noindent}\def\ov{\over}\def\p{\partial}
\def\smo{{\cal C}^\infty}\def\smoc{{\cal C}^\infty_c}\def\vp{\varphi}

\def\Gam{\Gamma}
\def\skp{\smallskip}

\newtheorem{theorem}{Theorem}[section]
\newtheorem{prop}{Proposition}[section]
\newtheorem{lemma}{Lemma}[section]
\newtheorem{cor}{Corollary}[section]

\title{Fundamental Solutions for the Tricomi Operator, II}
\author{J.~Barros-Neto\\
        Rutgers University, Hill Center, Piscataway, NJ 08954-8019\\
        \and
        I.~M.~Gelfand\\
        Rutgers University, Hill Center, Piscataway, NJ 08954-8019} 
\date{} 
\begin{document}  
\maketitle 
\begin{abstract}
In this paper we explicitly calculate fundamental solutions for the Tricomi operator, relative to an
arbitrary point in the plane, and show that all such fundamental solutions originate from the
hypergeometric function $F(1/6,1/6;1;\zeta)$ that is obtained when we look for homogeneous solutions to
the reduced hyperbolic Tricomi equation. 
\end{abstract}

\section{Introduction}\setcounter{equation}{0}\label{s1}
The Tricomi operator
\begin{equation}\label{eq1} 
{\cal T}=y{\p^2\over\p x^2}+{\p^2\over\p y^2},
\end{equation}
one of the simplest examples of a partial differential operator of {\em mixed type}, is: (i) {\em
elliptic} in the upper half plane ($y>0$);  {\em parabolic} along the $x$-axis ($y=0$); and (iii) {\em
hyperbolic} in the  lower half plane ($y<0$).

Our aim is to obtain explicit solutions in the sense of distributions or generalized functions of the equation
\begin{equation}\label{eq1c}
{\cal T}E=\delta(x-x_0,y-y_0),
\end{equation}
where $\delta(x-x_0,y-y_0)$ is the Dirac function at $(x_0,y_0)$ an arbitrary point in the plane. A solution 
$E$ of (\ref{eq1c}) is said to be a {\em fundamental solution relative to the point $(x_0,y_0)$.} 

In a previous paper \cite{bg} we considered the case when $x_0=y_0=0$ and proved the existence of two
remarkable fundamental solutions that clearly reflect the fact that the operator changes type across the
$x$-axis. In this paper we study the general case when $x_0=y_0\neq0$ and compare our results to those of
\cite{bg}. 

For sake of completeness and in order to make the reading of this paper independent of that of
\cite{bg}, we briefly review the contents of that paper. It is known that the equation $9x^2+4y^3=0$
defines the two characteristics for the Tricomi operator that emanate from the origin. These
characteristics, which are tangent to the $y$-axis at the origin, divide the plane in two disjoint regions
$D_+$ and $D_-$ (Figure 1) defined as follows:
$$
D_+=\{(x,y)\in{\mathbb R}^2 : 9x^2+4y^3>0\},
$$
the region ``outside'' the characteristics, and
$$
D_-=\{(x,y)\in{\mathbb R}^2 : 9x^2+4y^3<0\},
$$
the region ``inside'' the characteristics.
\begin{center}
\epsfig{file=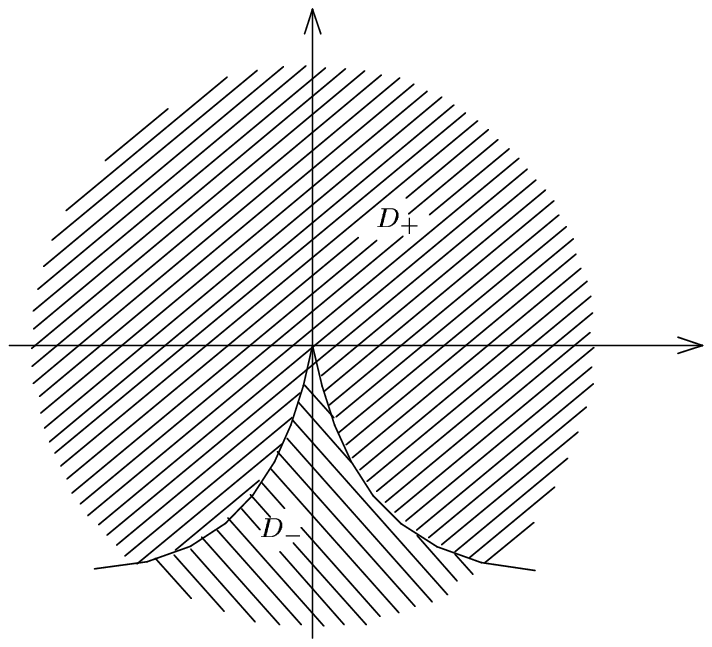, height=6cm, width=4cm}
\begin{center}
Figure 1
\end{center}
\end{center}

The first fundamental solution is defined by

\begin{equation}\label{eq1a}
\hskip .5cm F_+(x,y) =\left\{\begin{array}{ll}
         C_+(9x^2+4y^3)^{-1/6} & \mbox{in $D_+$}\\
\\
          \ \ \  0 & \mbox{elsewhere}
\end{array} \right .
\end{equation}
with
\begin{equation}\label{eq1a'}
C_+=-\ds{1\ov2^{1/3}\cdot3^{1/2}}F({1\ov6},{1\ov6};1;1),
\end{equation}
and the second one is defined by
\begin{equation}\label{eq1b}
F_-(x,y) = \left\{\begin{array}{ll}
         C_-|9x^2+4y^3|^{-1/6} & \mbox{in $D_-$}\\
\\
                          0 & \mbox{elsewhere}
\end{array} \right .
\end{equation}
with
\begin{equation}\label{eq1b'}
C_-=\ds{1\ov2^{1/3}}F({1\ov6},{1\ov6};1;1).
\end{equation}
In the expressions of both $C_+$ and $C_-,$ the constant
$$
F({1\ov6},{1\ov6};1;1)={\Gam(2/3)\ov\Gam^2(5/6)}
$$
is the value of the hypergeometric function $F(1/6,1/6;1;\zeta)$ at $\zeta=1.$

We remark that in \cite{bg} the constants $C_+$ and $C_-$ were respectively denoted by
$$
-{\Gam(1/6)\ov3\cdot2^{2/3}\pi^{1/2}\Gam(2/3)}\qquad\mbox{and}\qquad 
          {3\Gam(4/3)\ov2^{2/3}\pi^{1/2}\Gam(5/6)}.
$$
At that time, we were unaware of the role played by the hypergeometric function $F(1/6,1/6;1;\zeta)$ in
the study of fundamental solutions for the Tricomi operator. 

We also observe that according to formula (\ref{eq1b}), a perturbation at the origin spreads out to the entire
region $D_-,$ inside the characteristics. On the other hand, according to formula (\ref{eq1a}), the same
perturbation spreads out to the whole elliptic region ($y>0$) and also to the hyperbolic region ($y<0$)
outside the characteristics.  Since the Tricomi operator is invariant under translations along the
$x$-axis, the same phenomena take place for fundamental solutions that correspond to the Dirac measure
$\delta$ concentrated at an arbitrary point $(a,0)$ on the $x$-axis.

It is well know that the Tricomi operator describes the transition from subsonic flow (elliptic region) to
supersonic flow (hyperbolic region). In the extensive literature on the Tricomi operator one finds in the
works of several authors, among others, Agmon \cite{ag}, Friedrichs \cite{fri}, Gelfand [unpublished],
Germain and Bader \cite{gb1}, Landau and Lifshitz \cite{ll}, Leray \cite{ler2}, and Morawetz \cite{mo}, a
large number of examples of solutions to different problems that show the interaction between these two
regions. It seems that the two distributions $F_+(x,y)$ and $F_-(x,y)$ are the simplest of such examples.

We now return to equation (\ref{eq1c}). In view of the invariance of the Tricomi operator
under translations parallel to the $x$-axis,  the problem of solving that equation is equivalent to that of
solving 
\begin{equation}\label{eq1c'}
{\cal T}E=\delta(x,y-b),
\end{equation}
where $b$ is an arbitrary real number and $\delta(x,y-b)$ denotes the Dirac measure concentrated at the
point $(0,b).$ 

Consider the case $b<0$, that is, the point $(0,b)$ is situated in the hyperbolic region. We introduce
the {\em characteristic coordinates} 
\begin{equation}\label{ch0}
\ell=3x+2(-y)^{3/2}\quad\mbox{and}\quad m=3x-2(-y)^{3/2},
\end{equation}
set $\ell_0=2(-b)^{3/2},$ and to let $a>0$ be such that $3a=2(-b)^{3/2}.$ Note that $(\ell_0,-\ell_0)$
represents the point $(0,b)$ in characteristic coordinates. As shown 
\begin{center}
\epsfig{file=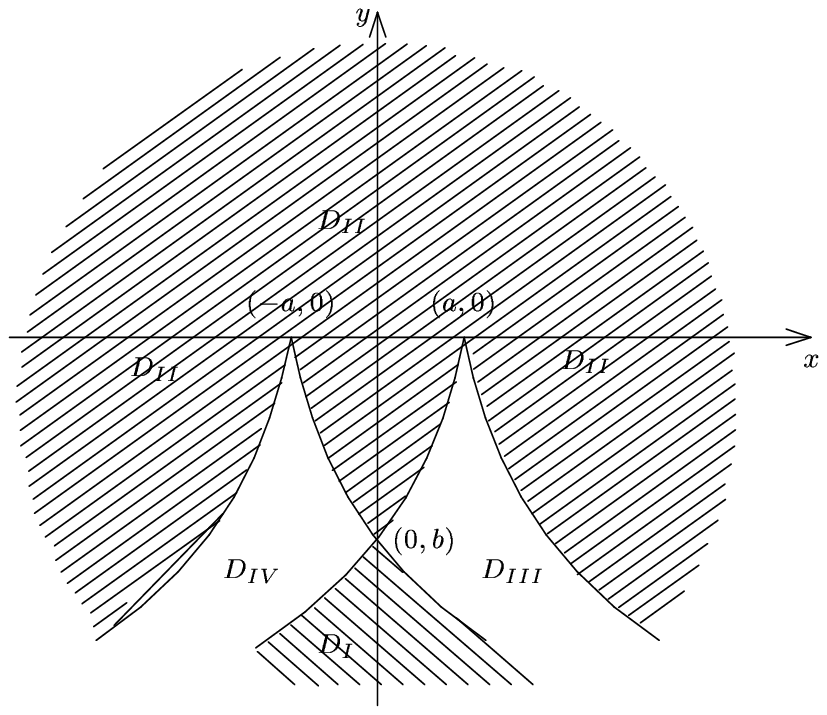, height=8cm, width=6cm}
\begin{center}
Figure 2
\end{center}
\end{center}
in Figure 2, two characteristics with
equations
\begin{equation}\label{ch1}
3(x-a)+2(-y)^{3/2}=0\quad\mbox{and}\quad 3(x+a)-2(-y)^{3/2}=0
\end{equation}
pass through the point $(0,b).$ In characteristic coordinates, these are the two half-lines
\begin{equation}\label{ch1a}
\ell=\ell_0,\ -\infty<m\leq\ell_0,\quad\mbox{and}\quad m=-\ell_0,\ -\ell_0\leq\ell<+\infty
\end{equation}
that originate from the point $(\ell_0,-\ell_0).$

The characteristics (\ref{ch1}) meet the $x$-axis, respectively, at the points $(a,0)$ and $(-a,0).$ Two new
characteristics originate from these two points, namely,
\begin{equation}\label{ch2}
3(x-a)-2(-y)^{3/2}=0\quad\mbox{and}\quad3(x+a)+2(-y)^{3/2}=0,
\end{equation}
which we call {\em reflected characteristics.} The corresponding equations in characteristic coordinates
are
\begin{equation}\label{ch2a}
m=\ell_0,\ \ell_0\leq\ell<+\infty\quad\mbox{and}\quad\ell=-\ell_0,\ -\infty<m\leq-\ell_0.
\end{equation}

The two characteristics through the point $(0,b)$ and the reflected characteristics divide the plane into
four disjoint regions denoted by $D_I,$ $D_{II},$ $D_{III},$ and $D_{IV},$ and illustrated in Figure 2. 

In this paper, we show existence of four fundamental solutions each supported by the closure of the
corresponding region. Of these solutions, only two have physical meaning: the one defined in the region $D_I$
and the other defined in the region $D_{II}.$ This situation is similar to what happens in the case of the
wave  operator in two dimensions where the two relevant fundamental solutions are the ones supported by
the forward and backward light-cones. As we will see, it is the hypergeometric function $F(1/6,1/6;1;\zeta)$ 
that plays a crucial role in defining these four fundamental solutions.

The plan of this paper is as follows.  We prove in Section \ref{s2} that the function
\begin{equation}\label{eq6}
\hskip 1cm E(\ell,m;\ell_0,-\ell_0)=(\ell+\ell_0)^{-1/6}(\ell_0-m)^{-1/6}
     F({1\ov6},{1\ov6};1;{(\ell-\ell_0)(m+\ell_0)\ov(\ell+\ell_0)(m-\ell_0)})
\end{equation}
is a solution of ${\cal T}_hu=0,$ where ${\cal T}_h$ denotes the {\em reduced hyperbolic} Tricomi equation 
(\ref{eq3}). 

After replacing into (\ref{eq6}), $\ell$ and $m$ by their expressions in (\ref{ch0}),  we obtain the
following function of $x$ and $y$:
\begin{eqnarray}\label{eq6'}
\hskip 1cm E(x,y;0,b)&=&e^{i\pi/6}[9(x^2-a^2)+4y^3-12a(-y)^{3/2}]^{-1/6}\times\\
& &\phantom{a}\times
F\left({1\ov6},{1\ov6};1;{9(x^2-a^2)+4y^3+12a(-y)^{3/2}\ov9(x^2-a^2)+4y^3-12a(-y)^{3/2}}\right).
\nonumber
\end{eqnarray}

This is the function that generates the four fundamental solutions relative to the point $(0,b).$
In Section \ref{s3} we define the distribution $E_I(x,y;0,b)$ as the restriction of $E(x,y;0,b),$ after multiplication by a
suitable constant, to the region $D_I$ and show that $E_I$ is a fundamental solution of the Tricomi operator. The
region $D_I$ is entirely contained in the hyperbolic region where it is natural to use characteristic coordinates.       
Theorem \ref{th2} is then proved, via integration by parts and by using the results of Proposition \ref{pr1}. As a 
consequence of it we show that, as $(0,b)$ tends to $(0,0),$ $E_I(x,y;0,b)$ tends, in the sense of distributions, to
$F_-(x,y),$ the fundamental solution defined by (\ref{eq1b}).   

The fundamental solutions supported by the closure of the regions $D_{III}$ and $D_{IV}$ are defined and studied
in Section \ref{s4}. In both regions, it is necessary to take into account the singularities of $E(x,y;0,b)$ along the
reflected characteristics. The results of Proposition \ref{pr2}, that describe the asymptotic behavior of the
hypergeometric function $F(1/6,1/6;1;\zeta)$ as $|\zeta|\to+\infty,$ are then needed in the proof of Theorem
\ref{th3}. We note that as $b\to 0,$ both fundamental solutions $E_{III}$ and $E_{IV}$ tend to 0.

In Section \ref{s5}, we define $E_{II}(x,y;0,b)$ as the restriction of $E(x,y;0,b),$ after multiplication by a suitable
constant, to $D_{II}$ and show that ${\cal T}E_{II}=\delta(0,b).$ Since $E_{II}$ is complex valued and the Tricomi
operator has real coefficients, both its complex conjugate and real part are also fundamental solutions relative to
the point $(0,b).$ The imaginary part of $E_{II}$ is then a solution of the homogeneous equation ${\cal T}u=0$ and we
will call   it a {\em Tricomi harmonic function}.

Contrary to what was proved in Corollary \ref{cor1}, it is not true that the fundamental solution $F_+(x,y)$
defined by (\ref{eq1a}) is the limit of $E_{II}(x,y;0,b),$ or its real part, as $(0,b)\to(0,0).$ It is necessary to take a suitable
linear combination of $E_{II}$ and its complex conjugate in order to achieve this result.

As a final remark, we mention that in his paper \cite{ler2}, J.~Leray described a general method, based upon the theory
of analytic functions of several complex variables to find fundamental solutions for a class of hyperbolic linear differential
operator with analytic coefficients. In particular, he showed how his method could be used to obtain, in the hyperbolic
region, a fundamental solution for the Tricomi operator relative to a point $(0,b).$ He also produced an explicit formula 
for the fundamental solution in terms of the hypergeometric function $F(1/6,1/6;1;\zeta)$. Our method is simpler, 
more direct, and gives us global fundamental solutions that clearly reflect the change of type of the Tricomi operator
across the $x$-axis.

We would like to thank Abbas Bahri, Fernando Cardoso and Vladimir Retakh for several helpful discussions.

\section{A Special Solution to ${\cal T}_hu=0$}\label{s2}\setcounter{equation}{0} 

In characteristic coordinates (\ref{ch0}) the Tricomi operator $\cal T$ becomes  
\begin{equation}\label{eq3}
{\cal T}_h = {\p^2\over\p\ell\p 
m}-{1/6\over{\ell-m}}({\p\over\p\ell}-{\p\over\p m}),
\end{equation}
and we call ${\cal T}_h$ the {\em reduced hyperbolic} form of $\cal T.$

We now look for homogeneous solutions of the equation ${\cal T}_hu=0.$ Every homogeneous  function of $\ell$ and
$m$ of degree $\lambda,$ a complex number, can be written as
$$
u(\ell,m)=\ell^\lambda\phi(t),
$$
where $\phi$ is a function of a single variable $t=m/\ell.$ Direct substitution into (\ref{eq3}) shows that
$\phi(t)$ must be a solution of the hypergeometric equation
$$
t(1-t)\phi''(t)+[({5\over 6}-\lambda)-({7\over6}-\lambda)t]\phi'(t) + 
{\lambda\over 6}\phi(t)=0.
$$
As a solution of this equation, we choose the following hypergeometric function
$F(-\lambda,1/6;5/6-\lambda;t)$ extended, by analytic continuation, to the whole complex plane
$\Bbb C$ minus the cut $[1,\infty).$ Since we are looking for fundamental solutions to the Tricomi operator we take
for $\lambda$ the value $-1/6,$ as previously indicated in our joint paper \cite{bg}. Thus
$$
u(\ell,m)=\ell^{-1/6}F({1\ov 6},{1\ov 6};1;{m\ov\ell})
$$
is a solution of ${\cal T}_hu=0.$

Let now $(\ell_0,m_0)$ be an arbitrary point in ${\Bbb R}^2$ and consider 
the change of variables
$$
\ell \to {\ell-m_0\ov \ell-\ell_0}\qquad m \to {m-m_0\ov m-\ell_0}.
$$
After unenlightening calculations, one can show that the function
\begin{equation}\label{eq5b} 
\hskip 1cm E(\ell,m;\ell_0,m_0)\!=\!(\ell-m_0)^{-1/6}(\ell_0-m)^{-1/6}
             F({1\ov 6},{1\ov 6};1;\!{(\ell-\ell_0)(m-m_0)\ov(\ell-m_0)(m-\ell_0)})
\end{equation}
is also a solution of equation (\ref{eq3}). This is the special solution to ${\cal T}_hu=0$ that we are looking for. 

Consider now the adjoint to equation (\ref{eq3}):
\begin{equation}\label{eq5c}
{\cal T}^*_h v = {\p^2v\over\p\ell\p m}+{1/6\over{l-m}}({\p 
v\over\p\ell}-{\p v\over\p m})
  -{1/3\ov(\ell-m)^2}v=0.
\end{equation}
One can see that if $v$ is a solution of equation (\ref{eq5c}), then $u=(\ell-m)^{-1/3}v$ is a solution of
equation (\ref{eq3}).  Since $E(\ell,m;\ell_0,m_0)$ is a solution of (\ref{eq6}), it follows that the function
\begin{equation}\label{eq5d}
R(\ell,m; \ell_0,m_0)=(\ell-m)^{1/3}E(\ell,m;\ell_0,m_0)
\end{equation}
is a solution of the adjoint equation (\ref{eq5c}). 

\begin{prop}\label{pr1}
$R(\ell,m;\ell_0,m_0)$ is the unique solution of ${\cal T}^*_hv=0$ that 
satisfies the following
conditions:\par
(i) $R_\ell=\ds{1/6\over{\ell-m}}R$ along the line $m=m_0,$\par
(ii) $R_m=\ds{-1/6\over{\ell-m}}R$ along the line $\ell=\ell_0$;\par
(iii) $R(\ell_0,m_0;\ell_0,m_0)=1.$
\end{prop}

\ni{\bf Proof.} Clearly, conditions $(i),\ (ii), \mbox{and} \ (iii)$ 
imply
uniqueness for $R.$ As $m=m_0,$ the argument of the hypergeometric 
function in (\ref{eq5b}) 
equals zero and so $F(1/6,1/6;1;0)=1.$ Thus, along the line $m=m_0,$ we 
have that
\begin{equation}\label{eq5e}
R(\ell,m_0;\ell_0,m_0)={\left(\ell-m_0\ov\ell_0-m_0\right)^{1/6}}
  =e^{\ds\int_{\ell_0}^\ell a(t)\,dt},
\end{equation}
where $a(t)=\ds{1/6\ov t-m_0}.$ Therefore $R_\ell=\ds{1/6\ov \ell-m}R$ 
along $m=m_0.$ In the
same manner, one can see that $R_m=\ds{-1/6\ov \ell-m}R$ along the line 
$\ell=\ell_0.$ Finally, it is clear that $R(\ell_0,m_0;\ell_0,m_0)=1.$
\hfill$\Box$ 

\medskip
{\bf Remark 1.} The function $R(\ell,m;\ell_0,m_0)$ is the Riemann function of the operator ${\cal
T}_h$ relative to the point $(\ell_0,m_0),$ (see \cite{ch}). 

\medskip
If we go back to the notations in the Introduction, with $(0,b)$ and $b<0,$ and let $\ell_0=2(-b)^{3/2},$ then
$m_0=-\ell_0$ and we obtain from equation (\ref{eq5b}) the equation (\ref{eq6}) that, for further reference, 
we write as follows:
\begin{equation}\label{ch3} 
\hskip .7cm 
E(\ell,m;\ell_0,-\ell_0)=(\ell+\ell_0)^{-1/6}(\ell_0-m)^{-1/6}
             F({1\ov 6},{1\ov 6};1;\zeta)
\end{equation}
where
\begin{equation}\label{ch3a}
\zeta={(\ell-\ell_0)(m+\ell_0)\ov(\ell+\ell_0)(m-\ell_0)}.
\end{equation}

Similarly, after replacing $\ell$ and $m$ by their expressions in (\ref{ch0}) we obtain the function (\ref{eq6'}) that we
rewrite as 
\begin{equation}\label{ch4}
\hskip 1cm E(x,y;0,b)=e^{i\pi/6}[9(x^2-a^2)+4y^3-12a(-y)^{3/2}]^{-1/6}\times
F\left({1\ov6},{1\ov6};1;\zeta\right)
\nonumber
\end{equation}
where
\begin{equation}\label{ch4a}
\zeta={9(x^2-a^2)+4y^3+12a(-y)^{3/2}\ov9(x^2-a^2)+4y^3-12a(-y)^{3/2}}.
\end{equation}

We wish to analyze this function. Denote by $r_a$ the reflected characteristic $3(x-a)-2(-y)^{3/2}=0$ at $(a,0)$ and
by $r_{-a}$ the reflected characteristic $3(x+a)+2(-y)^{3/2}=0$ at $(-a,0).$ 
\begin{prop}\label{pr1a}
$E(x,y;0,b)$ is a real analytic function defined in ${\mathbb R}^2\backslash(r_a\cup r_{-a}).$ 
The singularities of $E(x,y;0,b)$ occur along the reflected characteristics $r_a$ and $r_{-a}.$
\end{prop}

\ni{\bf Proof.} 1. Let
$$
z=[3(x+a)+2(-y)^{3/2}][3(x-a)-2(-y)^{3/2}]=9(x^2-a^2)+4y^3-12a(-y)^{3/2}
$$
be the denominator of $\zeta$ in the expression (\ref{ch4a}). For $y\leq0,$ $z=0$ on $r_a\cup r_{-a},$ it is
{\em negative} in the region inside these characteristics, and it is {\em positive} in the region outside them. For
$y>0,$ $z$ is a complex number. 

It is possible to choose the argument of $z$ so that it varies continuously in the region 
${\mathbb R}^2\backslash(r_a\cup r_{-a}).$ Indeed, define  
$\phi:{\mathbb R}^2\backslash(r_a\cup r_{-a})\to S_1,$ the unit circumference, by
$$
\phi(x,y)={9(x^2-a^2)+4y^3-12a(-y)^{3/2}\ov\rho},
$$
where 
$$\rho=|9(x^2-a^2)+4y^3-12a(-y)^{3/2}|.
$$ 
Since ${\mathbb R}^2\backslash(r_a\cup r_{-a})$ is contractible,  $\phi$ lifts to ${\mathbb R},$ that is, there 
exists a continuous function $\theta(x,y)$ defined on ${\mathbb R}^2\backslash(r_a\cup r_{-a})$ so that
$$
9(x^2-a^2)+4y^3-12a(-y)^{3/2}=\rho e^{-i\theta(x,y)}.
$$
We may take $\theta(x,y)$ to be equal zero outside the reflected characteristis and equal $\pi$ inside them.
It follows that
$$
z^{-1/6}=[9(x^2-a^2)+4y^3-12a(-y)^{3/2}]^{-1/6}
$$
is well defined and real analytic in ${\mathbb R}^2\backslash(r_a\cup r_{-a}).$ As for
$$
F\left({1\ov6},{1\ov6};1;
     {9(x^2-a^2)+4y^3+12a(-y)^{3/2}\ov 9(x^2-a^2)+4y^3-12a(-y)^{3/2}}\right)
$$
it is also real analytic in the same region. \hfill $\Box$
\medskip

The following remaks are in order:

1) For $y=0$ we have
$$
E(x,0;0,b)=F({1\ov6},{1\ov6}, 1;1)|x^2-a^2|^{-1/6},\quad \forall x\in(-a,a)
$$
and
$$
E(x,0;0,b)=e^{i\pi/6}F({1\ov6},{1\ov6}, 1;1)(x^2-a^2)^{-1/6},\quad \forall x\notin(-a,a).
$$

2) For $y>0$ we have
$$
\zeta={\bar z\ov z}=\rho e^{i2\theta(x,y)},\quad\mbox{so}\quad |\zeta|=1.
$$
It is known (see \cite{erd}) that the hypergeometric function $F(1/6,1/6;1;\zeta)$ is then given by its  
hypergeometric series and this series is absolutely convergent in the closed disk $|\zeta|\leq 1.$ 

3) Along the characteristics
$$
3(x-a)+2(-y)^{3/2}=0 \quad\mbox{and}\quad  3(x+a)-2(-y)^{3/2}=0
$$
through the point $(0,b)$ the function $E(x,y;0,b)$ is equal to
$$
E=2^{-2/3}(by)^{-1/4}
$$
and so it has a singularity of order $-1/4$ at $y=0.$

4) Along both reflected characteristics $r_{-a}$ and $r_{a}$, the function $E(x,y;0,b)$ has a logarithmic singularity.
This follows from Proposition \ref{pr2}, in  Section \ref{s4}, that describes the asymptotic behavior of the 
hypergeometric function $F(1/6,1/6;1;\zeta)$ as $|\zeta|\to+\infty.$ 

\section{The Fundamental Solution in the region $D_I$}\label{s3}\setcounter{equation}{0} 

Let $E(x,y;0,b)$ be the function defined by the expressions (\ref{ch4}) and (\ref{ch4a}) and define 
\begin{equation}\label{eq8}
E_I(x,y;0,b) = \left\{\begin{array}{ll}
         {\ds{1\over2^{1/3}}}E(x,y;0,b) & \mbox{in $D_I$}\\
  \\
                          0 & \mbox{elsewhere.}
\end{array} \right .
\end{equation}
Since $E(x,y;0,b)$ is $\smo$ in $D_I$ and bounded on the boundary of $D_I,$ it follows that $E_I(x,y;0,b)$ is
a locally integrable function and defines a distribution whose support is the closure of $D_I.$ 

\begin{theorem}\label{th2} 
$E_I(x,y;0,b)$ is a fundamental solution for the Tricomi operator $\cal T$ relative to the point $(0,b).$
\end{theorem}

\ni{\bf Proof.} We must show that
$$
\langle E_I,{\cal T}\varphi\rangle=\varphi(0,b),\quad 
         \forall\varphi\in\smoc({\Bbb R}^2_-).
$$
Since $E(x,y;0,b)$ is locally integrable in ${\Bbb R}^2,$ this is equivalent to showing that
\begin{equation}\label{int1}
\int\!\int_{{\Bbb R}^2_-}E_I(x,y,;0,b){\cal T}\varphi(x,y)\,dx\,dy=\varphi(0,b),
        \forall\phi\in\smoc({\Bbb R}^2_-).
\end{equation}

By introducing the characteristic coordinates $\ell$ and $m,$ noting that in the new variables the Tricomi
operator (\ref{eq1}) becomes
$$
y{\p^2\over\p x^2}+{\p^2\over\p
y^2}=-2^{2/3}3^2(\ell-m)^{2/3}\left[{\p^2\over\p\ell\p
m}-{1/6\over{l-m}}({\p\over\p\ell}-{\p\over\p m})\right],
$$
and that the Jacobian of the transformation is
$$
{\p(x,y)\over\p(\ell,m)}={1\over2^{1/3}3^2(\ell-m)^{1/3}},
$$
we can see that the integral in (\ref{int1}) is equal to
$$
\!-\!\int_{-\infty}^{-\ell_0}\!\!\int_{\ell_0}^{\infty}\!(\ell-m)^{1/3}
E(\ell,m;\ell_0,-\ell_0)
 (\phi_{\ell m}\!-\!{1/6\ov \ell-m}\phi_\ell\!+\!{1/6\ov 
\ell-m}\phi_m)\,d\ell\,dm,
$$
where $\phi(\ell,m)$ is $\varphi(x,y)$ in characteristic coordinates.  We denote the last integral by $I$ 
and write it as

\begin{equation}\label{int2}
\hskip .75cm I
=-\int_{-\infty}^{-\ell_0}\!\int_{\ell_0}^{\infty}R(\ell,m;\ell_0,-\ell_0)
   (\phi_{\ell m}-{1/6\ov \ell-m}\phi_\ell+{1/6\ov 
\ell-m}\phi_m)\,d\ell\,dm,  
\end{equation}
after recalling that 
$(\ell-m)^{1/3}E(\ell,m;\ell_0,-\ell_0)=R(\ell,m;\ell_0,-\ell_0).$  

\begin{center}
\epsfig{figure=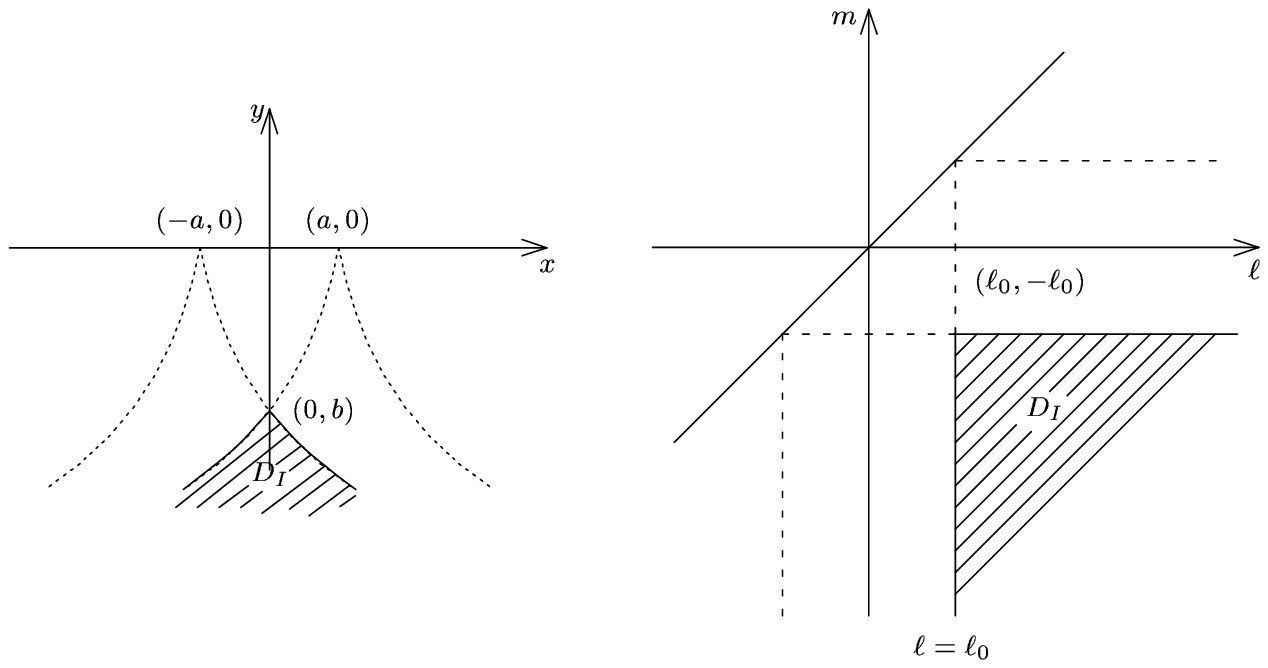, height=7.5cm, width=13.5cm}
\begin{center}
Figure 3
\end{center}
\end{center}

Next we perform several integrations by parts and take into account the properties of the function
$R(\ell,m;\ell_0,m_0),$ with $m_0=-\ell_0,$ as stated in Proposition \ref{pr1}. First consider the term
in  (\ref{int2}) that contains $\phi_{\ell m}.$ Integrating by parts first with respect to $m,$ say, and
then with respect to $\ell,$ one obtains
\begin{equation}\label{int3}
-\int_{-\infty}^{-\ell_0}\int_{\ell_0}^\infty R\phi_{\ell 
m}\,d\ell\,dm=\phi(\ell_0,-\ell_0)
\end{equation}
$$
+{1\ov6}\int_{\ell_0}^\infty{R(\ell,-\ell_0;\ell_0,-\ell_0)\ov\ell+\ell_0
}\phi(\ell,-\ell_0)\,d\ell   
+{1\ov6}\int_{-\infty}^{-\ell_0}{R(\ell_0,m;\ell_0,-\ell_0)\ov\ell_0-m}
\phi(\ell_0,m)\,dm
$$
$$
-\int_{-\infty}^{-\ell_0}\int_{\ell_0}^\infty R_{\ell m}\phi\,d\ell\,dm.
$$
 
Second integrate by parts, relative to $\ell,$ the term that contains 
$\phi_\ell$ in (\ref{int2}) and
obtain 
\begin{equation}\label{int4}
{1\ov6}\int_{-\infty}^{-\ell_0}\int_{\ell_0}^\infty{R\ov\ell-m}\phi_\ell\, d\ell\,dm=
\end{equation}
$$
-{1\ov6}\int_{-\infty}^{-\ell_0}{R(\ell_0,m;\ell_0,-\ell_0)\ov\ell_0-m}\phi(\ell_0,m)\,dm
-{1\ov6}\int_{-\infty}^{-\ell_0}\int_{\ell_0}^\infty{\p\ov\p\ell}
      \left({R\ov\ell-m}\right)\phi\,d\ell\,dm.
$$

Finally integrate by parts, relative to $m,$ the term that contains $\phi_m$ in (\ref{int2}) and get
\begin{equation}\label{int5}
-{1\ov6}\int_{-\infty}^{-\ell_0}\int_{\ell_0}^\infty{R\ov\ell-m}\phi_m\,d
\ell\,dm=
\end{equation}
$$
-{1\ov6}\int_{\ell_0}^{\infty}{R(\ell,-\ell_0;\ell_0,-\ell_0)\ov\ell+\ell
_0}\phi(\ell,-\ell_0)\,d\ell
+{1\ov6}\int_{-\infty}^{-\ell_0}\int_{\ell_0}^\infty{\p\ov\p m}
      \left({R\ov\ell-m}\right)\phi\,d\ell\,dm.
$$
By adding (\ref{int3}), (\ref{int4}), and (\ref{int5}) we obtain that 
$I=\phi(\ell_0,-\ell_0)$ and 
this completes the proof. \hfill$\Box$
\medskip

As a consequence of Theorem \ref{th2} we obtain the following result.

\begin{cor}\label{cor1}
As $(0,b)\to(0,0),$ the fundamental solution (\ref{eq8}) converges in the 
sense of distributions to $F_-(x,y)$ the fundamental solution defined by (\ref{eq1b}) and (\ref{eq1b'}).
\end{cor}

\ni{\bf Proof.} First note that whenever $\mbox{Re}(c-a-b)>0,$ the 
value $F(a,b;c;1)$ is finite and we have
$$
F(a,b;c;1)={\Gamma(c)\Gamma(c-a-b)\ov\Gamma(c-a)\Gamma(c-b)}.
$$ 
Thus $F(1/6,1/6;1;1)$ is well defined. As $(0,b)\to(0,0),$ then $(\ell_0,-\ell_0)\to(0,0),$ and so
the limit of $E(\ell,m;\ell_0,-\ell_0)$ defined by (\ref{ch3}) and (\ref{ch3a}) is
$$
(-\ell m)^{-1/6}F({1\ov6},{1\ov6};1;1).
$$
On the other hand, at the limit the region $D_I$ coincides with $D_-,$ the region inside the two
characteristics that originate from the origin, where $\ell m=9x^2+4y^3<0.$ \hfill$\Box$
\medskip

{\bf Remarks.}  1) In our joint paper \cite{bg}, the multiplicative constant appearing in the definition of 
$F_-(x,y)$ was given by 
$$
C_-={3\Gamma(4/3)\ov2^{2/3}\pi^{1/2}\Gamma(5/6)}.
$$
It is a matter of verification that this constant coincides with the one given by formula (\ref{eq1b'}).
 
2) Theorem \ref{th2} could have been proved in a different way, namely, by using the Green's formula
(\ref{gf}) and by integrating along a suitable contour, as we did  in our joint paper \cite{bg}. This method
which will be explained in the following section, is more adequate in proving existence of the 
fundamental solutions supported by the closure of the regions $D_{II},$ $D_{III},$ and $D_{IV}.$ 

\section{Fundamental solutions in the regions $D_{III}$ and  $D_{IV}$}\label{s4}\setcounter{equation}{0}  

For reasons of symmetry it is enough to consider only the region $D_{III}.$ This is the region bounded by part of the
characteristic $3(x+a)-2(-y)^{3/2}=0,$ part of the characteristic $3(x-a)+2(-y)^{3/2}=0$ and the reflected 
characteristic $3(x-a)-2(-y)^{3/2}=0$ (See Figure 2). 

As we already pointed out, the function $E(x,y;0,b)$  is singular along the reflected characteristic 
$3(x-a)-2(-y)^{3/2}=0.$ The nature of the singularity depends on the behavior of the hypergeometric function
$F(1/6,1/6;1;\zeta),$ where $\zeta=(\ell-\ell_0)(m+\ell_0)/(\ell+\ell_0)(m-\ell_0),$ along the
characteristic $m=\ell_0.$ This is revealed by the following result on the analytic continuation of the
hypergeometric series $F(a,a;c;\zeta)$ (see Erd\'ely \cite{erd}). The same result also gives
us the asymptotic behavior of the analytic extension as $|\zeta|\to\infty.$

\begin{prop}\label{pr2}
For $\zeta\in {\Bbb C}$ with $|\!\arg(-\zeta)|<\pi$, we have
\begin{equation}\label{eq10}
F(a,a;c;\zeta)=(-\zeta)^{-a}[\log(-\zeta)u(\zeta)+v(\zeta)]
\end{equation}
where
\begin{equation}\label{eq10a}
u(\zeta)={\Gamma(c)\ov\Gamma(a)\Gamma(c-a)}
   \sum_{n=0}^\infty{(a)_n(1-c+a)_n\ov(n!)^2}\zeta^{-n},
\end{equation}
with $(a)_n=\Gamma(a+n)/\Gamma(a),$ and
\begin{equation}\label{eq10b}
v(\zeta)={\Gamma(c)\ov\Gamma(a)\Gamma(c-a)}
   \sum_{n=0}^\infty{(a)_n(1-c+a)_n\ov(n!)^2}\,h_n\,\zeta^{-n},
\end{equation}
with $h_n=2\Psi(1+n)-\Psi(a+n)-\Psi(c-a-n),$ 
$\Psi(\zeta)=\Gamma'(\zeta)/\Gamma(\zeta).$ Moreover,
both series $u(\zeta)$ and $v(\zeta)$ converge for $|\zeta|>1.$
\end{prop}

From this proposition it follows, as we will show in the proof of Theorem \ref{th3}, that the singularity of
$E(x,y;0,b)$ along the reflected characteristics is logarithmic. We may now define

\begin{equation}\label{eq11}
E_{III}(x,y;0,b) = \left\{\begin{array}{ll}
         {\ds-{1\over2^{1/3}}}E(x,y;0,b) & \mbox{in $D_{III}$}\\
  \\
                          0 & \mbox{elsewhere,}
\end{array} \right .
\end{equation}
a distribution whose support is the closure of $D_{III}.$ We then have the following result 
\begin{theorem}\label{th3} 
$E_{III}(x,y;0,b)$ is a fundamental solution for the Tricomi operator $\cal T$ relative to $(0,b).$
\end{theorem}

\ni{\bf Proof.} 1.\ We use the following Green's formula for the Tricomi operator (see\cite{bg}):
\begin{equation}\label{gf}
\hskip .75cm
\int\!\!\int_D(E{\cal T}\vp-\vp{\cal T}E)\,dx\,dy=
 \int_CE(y\vp_x\,dy-\vp_y\,dx)-\vp(yE_x\,dy-E_y\,dx)
\end{equation}
where $D$ is a bounded domain with smooth boundary $C.$ If ${\cal T}E=0$ on D, then last formula becomes
\begin{equation}\label{gf1}
\int\!\!\int_DE{\cal T}\vp\,dx\,dy=
 \int_CE(y\vp_x\,dy-\vp_y\,dx)-\vp(yE_x\,dy-E_y\,dx)
\end{equation}
and the contour integral may be used to evaluate the double integral. Now throughout this section the 
domain $D=D_{III}$ lies entirely in the hyperbolic region and so it is more convenient  to express the
countour integral on the right-hand side of (\ref{gf1}) in terms of characteristic coordinates.

One can check that, in these coordinates, we have 
$$
y\vp_x\,dy-\vp_y\,dx=A(\ell-m)^{1/3}\psi_\ell\,d\ell-A(\ell-m)^{1/3}\psi_
m\,dm,
$$
where
$$
\psi(\ell,m)=\vp({\ell+m\ov6},-({\ell-m\ov4})^{2/3})
$$
and $A=1/2^{2/3}.$ Similarly,
$$
yE_x\,dy-E_y\,dx=A(\ell-m)^{1/3}E_\ell\,d\ell-A(\ell-m)^{1/3}E_m\,dm.
$$
Thus the contour integral in (\ref{gf1}), denoted by $I_C,$ can be written as
\begin{equation}\label{ci}
\hskip 1cm I_C = A\int_C(\ell-m)^{1/3}(E\psi_\ell-\psi E_\ell)\,d\ell
   -  A\int_C(\ell-m)^{1/3}(E\psi_m-\psi E_m)\,dm.\nonumber
\end{equation}

\begin{center}
\epsfig{figure=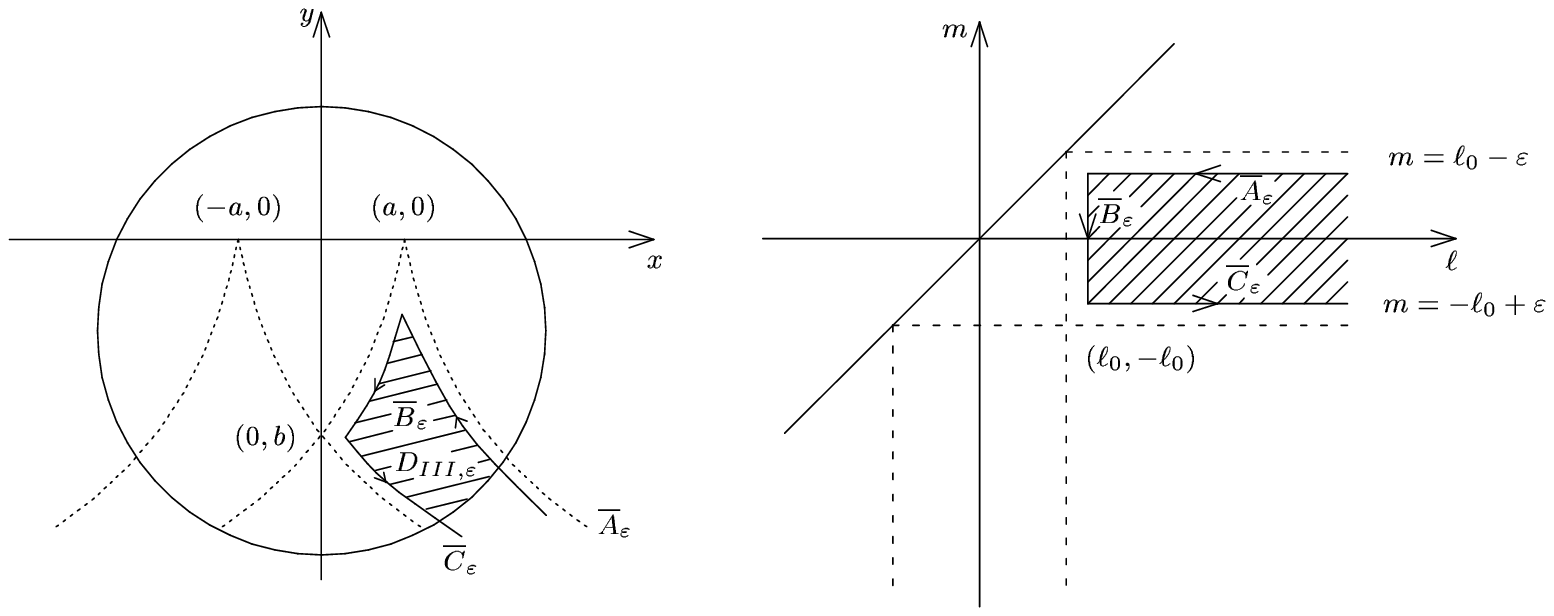, height=5.6cm, width=14cm}
\end{center}
\begin{center}
Figure 4
\end{center}

2.\ We must show that
\begin{equation}\label{eq12}
\hskip 0.5cm \langle E,{\cal T}\vp\rangle=
    \lim_{\epsilon\to0}\int\!\!\int_{D_\epsilon}E{\cal 
T}\vp\,dx\,dy=-2^{1/3}\vp(0,b),\quad
     \forall\vp\in\smoc({\Bbb R}^2),
\end{equation}
where $D_\epsilon$ is the intersection of an open disk that contains the 
support of $\vp$ and the region, defined in characteristic coordinates, as follows 
$$
D_{III,\epsilon}=\{(\ell,m):\ell>\ell_0+\epsilon,-\ell_0+\epsilon<m<\ell_
0-\epsilon\}.
$$
Since $\vp\equiv0$ on the boundary of the disk, it follows from Green's formula, that we can replace the
double integral in (\ref{eq12}) by an integral of the form (\ref{ci}), where the contour $C$ consists of
the oriented line segments $\bar A_\epsilon,$ $\bar B_\epsilon,$ and $\bar C_\epsilon$ shown in Figure 4.

We start by evaluating the simplest of these integrals, namely, the ones along  $\bar  B_\epsilon$ and 
$\bar C_\epsilon.$
\skp

3.\ {\em The integral along $\bar B_\epsilon.$} Denote by  $I_{\bar 
B_\epsilon}$ the integral along the line
segment, $\ell=\ell_0+\epsilon$ with $m$ varying from $\ell_0-\epsilon$ 
to $-\ell_0+\epsilon.$ Taking
into account the orientation of $\bar B_\epsilon,$ we obtain from formula 
(\ref{ci}) that
\begin{eqnarray}\label{eq13}
\hskip 1cm 
\lefteqn{I_{\bar B_\epsilon}= A\int_{-\ell_0+\epsilon}^{\ell_0-\epsilon}
    (\ell_0+\epsilon-m)^{1/3}E(\ell_0+\epsilon,m;\ell_0,-\ell_0)
             \psi_m(\ell_0+\epsilon,m)\,dm}\\
&& -A\int_{-\ell_0+\epsilon}^{\ell_0-\epsilon}
    (\ell_0+\epsilon-m)^{1/3}E_m(\ell_0+\epsilon,m;\ell_0,-\ell_0)
             \psi(\ell_0+\epsilon,m)\,dm,\nonumber
\end{eqnarray}
where $A=1/2^{2/3}.$

Integrating by parts the first integral in (\ref{eq13}), we get
\begin{equation}\label{eq14}
I_{\bar B_\epsilon} = A(\ell_0+\epsilon-m)^{1/3}E(\ell_0+\epsilon,m;
\ell_0,-\ell_0)\psi(\ell_0+\epsilon,m)\Bigr|^{\ell_0-\epsilon}_{-\ell_0+
 \epsilon} 
\end{equation}
$$
 -{A\ov3}\int_{-\ell_0+\epsilon}^{\ell_0-\epsilon}
    (\ell_0+\epsilon-m)^{-2/3}E(\ell_0+\epsilon,m;\ell_0,-\ell_0)
             \psi(\ell_0+\epsilon,m)\,dm
$$
$$
 -2A\int_{-\ell_0+\epsilon}^{\ell_0-\epsilon}
    (\ell_0+\epsilon-m)^{1/3}E_m(\ell_0+\epsilon,m;\ell_0,-\ell_0)
             \psi(\ell_0+\epsilon,m)\,dm 
$$
From equation (\ref{eq5d}) it follows that
$$
R_m=-{1\ov3}(\ell-m)^{-2/3}E (\ell,m;\ell_0,-\ell_0)+
   (\ell-m)^{1/3}E_m (\ell,m;\ell_0,-\ell_0).
$$
Substituting into (\ref{eq14}), we obtain
\begin{equation}\label{eq15}
I_{\bar B_\epsilon} = A(\ell_0+\epsilon-m)^{1/3}E(\ell_0+\epsilon,m;
\ell_0,-\ell_0)\psi(\ell_0+\epsilon,m)\Bigr|^{\ell_0-\epsilon}_{-\ell_0+
 \epsilon} 
\end{equation}
$$
-2A\int_{-\ell_0+\epsilon}^{\ell_0-\epsilon}[R_m(\ell_0+\epsilon,m;
\ell_0,-\ell_0) +
{1\ov6}{R(\ell_0+\epsilon,m;\ell_0,-\ell_0)\ov(\ell_0+\epsilon-m)}]
  \psi(\ell_0+\epsilon,m)\,dm.
$$
Now, as $\epsilon\to0,$ the last integral tends to zero because, by Proposition \ref{pr1},
$R_m=-R/6(\ell-m)$ along the line $\ell=\ell_0.$ On the other hand, the first term in (\ref{eq15}),
namely,  
\begin{eqnarray*}
\lefteqn{A(\ell_0+\epsilon-m)^{1/3}E(\ell_0+\epsilon,m;
\ell_0,-\ell_0)\psi(\ell_0+\epsilon,m)\Bigr|^{\ell_0-\epsilon}_{-\ell_0
+\epsilon}=}\\
& & \hskip 2cm 
A(2\epsilon)^{1/3}E(\ell_0+\epsilon,\ell_0-\epsilon;\ell_0,-\ell_0)
        \psi(\ell_0+\epsilon,\ell_0-\epsilon) -\\
& & \hskip 3cm 
A(2\ell_0)^{1/3}E(\ell_0+\epsilon,-\ell_0+\epsilon;\ell_0,-\ell_0)
        \psi(\ell_0+\epsilon,-\ell_0+\epsilon).
\end{eqnarray*}
tends to 
$-A(2\ell_0)^{1/3}E(\ell_0,-\ell_0;\ell_0,-\ell_0)\psi(\ell_0,-\ell_0),$  
as $\epsilon\to0.$
Therefore,
\begin{equation}\label{eq16}
\lim_{\epsilon\to0}I_{\bar 
B_\epsilon}={-1\ov2^{2/3}}\psi(\ell_0,-\ell_0).
\end{equation}
because $E(\ell_0,-\ell_0;\ell_0,-\ell_0)=(2\ell_0)^{-1/3}$ and 
$A=1/2^{2/3}.$
\skp

4.\ {\em The integral along $\bar C_\epsilon.$} Denote by $I_{\bar 
C_\epsilon}$ the integral along the line
$m=-\ell_0+\epsilon$ with $\ell_0+\epsilon<\ell<+\infty.$ From formula 
(\ref{ci}) we now obtain
\begin {eqnarray}\label{eq17}
\hskip 1cm 
\lefteqn{I_{\bar C_\epsilon}= A\int_{\ell_0+\epsilon}^{+\infty}
    (\ell+\ell_0-\epsilon)^{1/3}E(\ell,-\ell_0+\epsilon;\ell_0,-\ell_0)
             \psi_\ell(\ell,-\ell_0+\epsilon)\,d\ell}\\
&& -A\int_{\ell_0+\epsilon}^{+\infty}
(\ell+\ell_0-\epsilon)^{1/3}E_\ell(\ell,-\ell_0+\epsilon;\ell_0,-\ell_0)
             \psi(\ell,-\ell_0+\epsilon)\,d\ell,\nonumber
\end{eqnarray}
 As before, integrate by parts the first integral, use the formula
$$
R_\ell={1\ov3}(\ell-m)^{-2/3}E(\ell,m;\ell_0,-\ell_0) +
               (\ell-m)^{1/3}E_\ell(\ell,m;\ell_0,-\ell_0).
$$
to substitute for $E_\ell$ in the second integral, and rewrite (\ref{eq17}) as
$$
I_{\bar C_\epsilon} = 
-A(2\ell_0)^{1/3}E(\ell_0+\epsilon,-\ell_0,-\epsilon;\ell_0,-\ell_0)
\psi(\ell_0+\epsilon,-\ell_0+\epsilon)
$$
$$
-2A\int_{\ell_0+\epsilon}^{+\infty}[R_\ell(\ell,-\ell_0+\epsilon;\ell_0,-
\ell_0) -
 {1\ov6}{R(\ell,-\ell_0+\epsilon;\ell_0,-\ell_0)\ov(\ell+\ell_0-\epsilon)
}]
     \psi(\ell,-\ell_0-\epsilon)\,d\ell.
$$

As before, we obtain
\begin{equation}\label{eq19}
\lim_{\epsilon\to0}I_{\bar 
C_\epsilon}={-1\ov2^{2/3}}\psi(\ell_0,-\ell_0).
\end{equation}

5.\ {\em The integral along $\bar A_\epsilon.$} In this case, 
$m=\ell_0-\epsilon,$ $\ell$ varies from
$+\infty$ to $\ell_0+\epsilon,$ and the integral to be considered is 
\begin {eqnarray}\label{eq20}
\hskip 1cm 
\lefteqn{I_{\bar A_\epsilon}= -A\int_{\ell_0+\epsilon}^{+\infty}
    (\ell-\ell_0+\epsilon)^{1/3}E(\ell,\ell_0-\epsilon;\ell_0,-\ell_0)
             \psi_\ell(\ell,\ell_0-\epsilon)\,d\ell}\\
&& +A\int_{\ell_0+\epsilon}^{+\infty}
(\ell-\ell_0+\epsilon)^{1/3}E_\ell(\ell,\ell_0-\epsilon;\ell_0,-\ell_0)
             \psi(\ell,\ell_0-\epsilon)\,d\ell.\nonumber
\end{eqnarray}

Our aim is to show that

\begin{equation}\label{eq21}
\lim_{\epsilon\to0}I_{\bar A_\epsilon}=0.
\end{equation}
Once this is done, then by adding (\ref{eq16}), (\ref{eq19}), and 
(\ref{eq21}), we obtain (\ref{eq12}) and
the theorem will be proved. 

As we did in the previous two cases, we integrate by parts the first 
integral in (\ref{eq20}) and get
\begin{equation}\label{eq22}
I_{\bar 
A_\epsilon}=A(2\epsilon)^{1/3}E(\ell_0+\epsilon,\ell_0-\epsilon;\ell_0,-\ell_0)
    \psi(\ell_0+\epsilon,\ell_0-\epsilon)
\end{equation}
$$
+{A\ov3}\int_{\ell_0+\epsilon}^{+\infty}(\ell-\ell_0+\epsilon)^{-2/3}
E(\ell,\ell_0-\epsilon;\ell_0,-\ell_0)\psi(\ell,\ell_0-\epsilon)\,d\ell
$$
$$
+2A\int_{\ell_0+\epsilon}^{+\infty}(\ell-\ell_0+\epsilon)^{1/3}
E_\ell(\ell,\ell_0-\epsilon;\ell_0,-\ell_0)\psi(\ell,\ell_0-\epsilon)\,d\ell.
$$

In order to prove (\ref{eq21}) we now must take into account the asymptotic behavior of
both $E(\ell,m;\ell_0,-\ell_0)$ and $E_\ell(\ell,m;\ell_0,-\ell_0)$ for $m$ near $\ell_0$ and, more
specifically, the behavior of the hypergeometric function $F(1/6,1/6;1;\zeta)$ and its derivative
\begin{equation}\label{eq22a}
F'({1\ov6},{1\ov6};1;\zeta)={1\ov36}F({7\ov6},{7\ov6};2;\zeta),
\end{equation}
as $\zeta\to\infty.$

To see this, we start by expressing (\ref{eq22}) in terms of $F(1/6,1/6;1;\zeta)$ and its
derivative. In order to simplify notations we write from now on, 
$F(\zeta)=F(1/6,1/6;1;\zeta)$ and $G(\zeta)=F({7/6},{7/6};2;\zeta).$ In the case we are considering 
$m=\ell_0-\epsilon$ so, in view of formula (\ref{ch3a}), we set
$$
\zeta={(2\ell_0-\epsilon)(\ell-\ell_0)\ov-\epsilon(\ell+\ell_0)}={\sigma(
\epsilon)\ov-\epsilon}
\qquad\mbox{with}\qquad 
   \sigma(\epsilon)={(2\ell_0-\epsilon)(\ell-\ell_0)\ov(\ell+\ell_0)}.
$$
Also, recalling formula (\ref{ch3}), we have
\begin{equation}\label{eq22b}
E(\ell,\ell_0-\epsilon;\ell_0,-\ell_0)=\epsilon^{-1/6}(\ell+\ell_0)^{-1/6}
        F({\sigma(\epsilon)\ov-\epsilon}),
\end{equation}
and so

\begin{equation}\label{eq22c}
E_\ell(\ell,\ell_0-\epsilon;\ell_0,-\ell_0)=
\end{equation}
$$
 -{\epsilon^{-1/6}\ov6}(\ell+\ell_0)^{-7/6}F({\sigma(\epsilon)\ov-\epsilon})
 -{\epsilon^{-7/6}\ell_0(2\ell_0-\epsilon)\ov18}(\ell+\ell_0)^{-13/6}
           G({\sigma(\epsilon)\ov-\epsilon}),
$$
in view of (\ref{eq22a}) and taking into account that
$$
{d\zeta\ov d\ell}={2\ell_0(2\ell_0-\epsilon)\ov-\epsilon(\ell+\ell_0)^2}.
$$
By substituting (\ref{eq22b}) and (\ref{eq22c}) into (\ref{eq22}) and 
combining the resulting integrals we
obtain
\begin{equation}\label{eq23}
I_{\bar A_\epsilon}=A(2\epsilon)^{1/3}(2\ell_0+\epsilon)^{-1/6}
       \epsilon^{-1/6}F({\epsilon-2\ell_0\ov2\ell_0+\epsilon})
    \psi(\ell_0+\epsilon,\ell_0-\epsilon)
\end{equation}
$$
+{A\ov3}(2\ell_0-\epsilon)^{5/6}\int_{\ell_0+\epsilon}^{+\infty}
    {(\ell-\ell_0+\epsilon)^{-2/3}\ov(\ell-\ell_0)^{1/6}(\ell+\ell_0)}
({\sigma(\epsilon)\ov\epsilon})^{1/6}F({\sigma(\epsilon)\ov-\epsilon})
    \psi(\ell,\ell_0-\epsilon)\,d\ell
$$
$$
-{A\ov9}\ell_0(2\ell_0-\epsilon)^{-1/6}\int_{\ell_0+\epsilon}^{+\infty}
   {(\ell-\ell_0+\epsilon)^{1/3}\ov(\ell-\ell_0)^{7/6}(\ell+\ell_0)}
({\sigma(\epsilon)\ov\epsilon})^{7/6}G({\sigma(\epsilon)\ov-\epsilon})
    \psi(\ell,\ell_0-\epsilon)\,d\ell.
$$

As $\epsilon\to 0,$ the first term on the right-hand side of (\ref{eq23}) 
clearly tends to zero. We must
prove that the difference of the two integrals in (\ref{eq23}) also tends 
to zero as $\epsilon\to 0.$ For
this, according to formulas (\ref{eq10}), (\ref{eq10a}), and 
(\ref{eq10b}), we have
$$
({\sigma(\epsilon)\ov\epsilon})^{1/6}F({\sigma(\epsilon)\ov-\epsilon})=
\log({\sigma(\epsilon)\ov\epsilon})u({\sigma(\epsilon)\ov-\epsilon})+
    v({\sigma(\epsilon)\ov-\epsilon})
$$
and
$$
({\sigma(\epsilon)\ov\epsilon})^{7/6}G({\sigma(\epsilon)\ov-\epsilon})=
\log({\sigma(\epsilon)\ov\epsilon})U({\sigma(\epsilon)\ov-\epsilon})+
  V({\sigma(\epsilon)\ov-\epsilon}).
$$
By using these two expressions we combine the two integrals in 
(\ref{eq23}) and write them as a sum 
$I_{\bar A_\epsilon}^{(1)}$ + $I_{\bar A_\epsilon}^{(2)}$ - 
$\log\epsilon\cdot I_{\bar A_\epsilon}^{(3)},$  
where
\begin{equation}\label{eq24}
I_{\bar A_\epsilon}^{(1)}=
\end{equation}
$$
{A\ov3}(2\ell_0-\epsilon)^{5/6}\int_{\ell_0+\epsilon}^{+\infty}
    {(\ell-\ell_0+\epsilon)^{-2/3}\ov(\ell-\ell_0)^{1/6}(\ell+\ell_0)}
    \log(\sigma(\epsilon))u({\sigma(\epsilon)\ov-\epsilon})
    \psi(\ell,\ell_0-\epsilon)\,d\ell
$$
$$
\hskip 0.5cm 
-{A\ov9}\ell_0(2\ell_0-\epsilon)^{-1/6}\int_{\ell_0+\epsilon}^{+\infty}
   {(\ell-\ell_0+\epsilon)^{1/3}\ov(\ell-\ell_0)^{7/6}(\ell+\ell_0)}
     \log(\sigma(\epsilon))U({\sigma(\epsilon)\ov-\epsilon})
    \psi(\ell,\ell_0-\epsilon)\,d\ell,
$$

\begin{equation}\label{eq25}
I_{\bar A_\epsilon}^{(2)}=
\end{equation}
$$
{A\ov3}(2\ell_0-\epsilon)^{5/6}\int_{\ell_0+\epsilon}^{+\infty}
    {(\ell-\ell_0+\epsilon)^{-2/3}\ov(\ell-\ell_0)^{1/6}(\ell+\ell_0)}
    v({\sigma(\epsilon)\ov-\epsilon})\psi(\ell,\ell_0-\epsilon)\,d\ell
$$
$$
\hskip 0.5cm 
-{A\ov9}\ell_0(2\ell_0-\epsilon)^{-1/6}\int_{\ell_0+\epsilon}^{+\infty}
   {(\ell-\ell_0+\epsilon)^{1/3}\ov(\ell-\ell_0)^{7/6}(\ell+\ell_0)}
      V({\sigma(\epsilon)\ov-\epsilon})\psi(\ell,\ell_0-\epsilon)\,d\ell,
$$
and

\begin{equation}\label{eq26}
I_{\bar A_\epsilon}^{(3)}=
\end{equation}
$$
{A\ov3}(2\ell_0-\epsilon)^{5/6}\int_{\ell_0+\epsilon}^{+\infty}
    {(\ell-\ell_0+\epsilon)^{-2/3}\ov(\ell-\ell_0)^{1/6}(\ell+\ell_0)}
    u({\sigma(\epsilon)\ov-\epsilon})\psi(\ell,\ell_0-\epsilon)\,d\ell
$$
$$
\hskip 0.5cm 
-{A\ov9}\ell_0(2\ell_0-\epsilon)^{-1/6}\int_{\ell_0+\epsilon}^{+\infty}
   {(\ell-\ell_0+\epsilon)^{1/3}\ov(\ell-\ell_0)^{7/6}(\ell+\ell_0)}
      U({\sigma(\epsilon)\ov-\epsilon})\psi(\ell,\ell_0-\epsilon)\,d\ell.
$$

\begin{lemma}\label{lm1}
Each of the integrals $I_{\bar A_\epsilon}^{(i)}, 1\leq i\leq 3,$ is  
$O(\epsilon^{1/6})$.
\end{lemma}

\ni{\bf Proof of lemma.} Consider first the integral $I_{\bar 
A_\epsilon}^{(1)}.$ Since the series in
(\ref{eq10a}) converge for large values of $|\zeta|,$ we may write that
$$
u({\sigma(\epsilon)\ov-\epsilon})=u_0+{\epsilon\ov\sigma(\epsilon)}
\tilde u({\sigma(\epsilon)\ov-\epsilon}) \qquad\mbox{and}\qquad 
U({\sigma(\epsilon)\ov-\epsilon})=U_0+{\epsilon\ov\sigma(\epsilon)}
\tilde U({\sigma(\epsilon)\ov-\epsilon}),
$$
where $\tilde u$ and $\tilde U$ are bounded functions for small values of 
$\epsilon.$ 

By substituting these two expressions into (\ref{eq24}), rewrite it as a 
sum of three terms:
\begin{equation}\label{eq27}
I_{\bar A_\epsilon}^{(1)}=
\end{equation}
$$
\{ {A\ov3}(2\ell_0-\epsilon)^{5/6}u_0\int_{\ell_0+\epsilon}^{+\infty}
    {(\ell-\ell_0+\epsilon)^{-2/3}\ov(\ell-\ell_0)^{1/6}(\ell+\ell_0)}
    \log(\sigma(\epsilon))\psi(\ell,\ell_0-\epsilon)\,d\ell
$$
$$
\hskip 0.5cm 
-{A\ov9}\ell_0(2\ell_0-\epsilon)^{-1/6}U_0\int_{\ell_0+\epsilon}^{+\infty
}
   {(\ell-\ell_0+\epsilon)^{1/3}\ov(\ell-\ell_0)^{7/6}(\ell+\ell_0)}
     \log(\sigma(\epsilon))\psi(\ell,\ell_0-\epsilon)\,d\ell\} 
$$
$$
+{A\ov3}(2\ell_0-\epsilon)^{5/6}\int_{\ell_0+\epsilon}^{+\infty}
    {(\ell-\ell_0+\epsilon)^{-2/3}\ov(\ell-\ell_0)^{1/6}(\ell+\ell_0)}
    \log(\sigma(\epsilon)){\epsilon\ov\sigma(\epsilon)}\tilde 
u({\sigma(\epsilon)\ov-\epsilon})
    \psi(\ell,\ell_0-\epsilon)\,d\ell
$$
$$
\hskip 0.5cm 
-{A\ov9}\ell_0(2\ell_0-\epsilon)^{-1/6}\int_{\ell_0+\epsilon}^{+\infty}
   {(\ell-\ell_0+\epsilon)^{1/3}\ov(\ell-\ell_0)^{7/6}(\ell+\ell_0)}
     \log(\sigma(\epsilon)){\epsilon\ov\sigma(\epsilon)}\tilde 
U({\sigma(\epsilon)\ov-\epsilon})
    \psi(\ell,\ell_0-\epsilon)\,d\ell.
$$
The first term (inside the brackets) in (\ref{eq27}) tends to $0$ as 
$\epsilon\to 0.$ Indeed, its limit as 
$\epsilon\to 0$ is
$$
[{A(2\ell_0)^{5/6}u_0\ov3}-{A\ell_0(2\ell_0)^{-1/6}U_0\ov9}]
   \int_{\ell_0}^\infty{(\ell-\ell_0)^{-5/6}\ov\ell+\ell_0}
      \log({2\ell_0(\ell-\ell_0)\ov\ell+\ell_0})\psi(\ell,\ell_0)\,d\ell.
$$
The integral on the right-hand side is finite. On the other hand, 
$$
u_0={1\ov\Gamma(1/6)\Gamma(5/6)} \quad\mbox{and}\quad
U_0={1\ov\Gamma(5/6)\Gamma(7/6)}={6\ov\Gamma(1/6)\Gamma(5/6)},
$$
are the constant terms of the series $u(\zeta)$  and $U(\zeta)$ given by 
(\ref{eq10a}).
It is a matter of verification that the quantity inside the brackets is 
equal zero.

Next, consider the second term in (\ref{eq27}). After the change of 
variables $\ell-\ell_0=t\epsilon,$ one
can see that the absolute value of that term is bounded above by
$$
C\epsilon^{1/6}\int_1^{+\infty}(t+1)^{-2/3}t^{-7/6}\log t\,dt,
$$
with $C$ a constant independent of $\epsilon.$ With a similar 
calculation, one can show that the third
term in (\ref{eq27}) is bounded above by
 $$
C\epsilon^{1/6}\int_1^{+\infty}(t+1)^{1/3}t^{-13/6}\log t\,dt,
$$
with $C$ another constant. Therefore, $I_{\bar 
A_\epsilon}^{(1)}=O(\epsilon^{1/6}).$

The expression (\ref{eq26}) for $I_{\bar A_\epsilon}^{(3)},$ similar to 
that of $I_{\bar A_\epsilon}^{(1)},$
is even simpler because of the absence of the factor 
$\log(\sigma(\epsilon)).$ Thus, with a similar proof
we also  conclude that $I_{\bar A_\epsilon}^{(3)}=O(\epsilon^{1/6}).$ 
Finally, the integral $I_{\bar A_\epsilon}^{(2)}$ given by (\ref{eq25}) is analogous to 
$I_{\bar A_\epsilon}^{(3)}$ with $u$ and $U$
replaced by $v$ and $V$ the power  series defined by (\ref{eq10b}). The 
proof of the lemma is then complete. \hfill$\Box$
\skp

Recalling the expression (\ref{eq24}) for $I_{\bar A_\epsilon}$ one sees that Lemma \ref{lm1} implies 
(\ref{eq21}) which completes the proof of the theorem. \hfill$\Box$

\section{Fundamental Solutions in the region $D_{II}$}\label{s5}\setcounter{equation}{0} 

In the same manner as we did in the previous sections, define the distribution
\begin{equation}\label{eq30}
E_{II}(x,y;0,b) = \left\{\begin{array}{ll}
         {\ds{1\over2^{1/3}}}E(x,y;0,b) & \mbox{in $D_{II}$}\\
  \\
                          0 & \mbox{elsewhere,}
\end{array} \right .
\end{equation}
whose support is the closure of the region $D_{II}.$ Our aim is to prove the following result.

\begin{theorem}\label{th5}
$E_{II}(x,y;0,b)$ is a fundamental solution for the Tricomi operator relative to the point $(0,b).$
\end{theorem}

\ni{\bf Proof.} The proof is, with few modifications, analogous to that of Theorem
\ref{th3}.  As before, it suffices to show that 
\begin{equation}\label{eq31}
\hskip 0.5cm \langle E,{\cal T}\vp\rangle=
    \lim_{\epsilon\to0}\int\!\!\int_{D_\epsilon}E{\cal 
T}\vp\,dx\,dy=2^{1/3}\vp(0,b),\quad
     \forall\vp\in\smoc({\Bbb R}^2),
\end{equation}
where the domain of integration $D_\epsilon$ is defined as follows. Let $D$ be an open 

\begin{center}
\centering\epsfig{figure=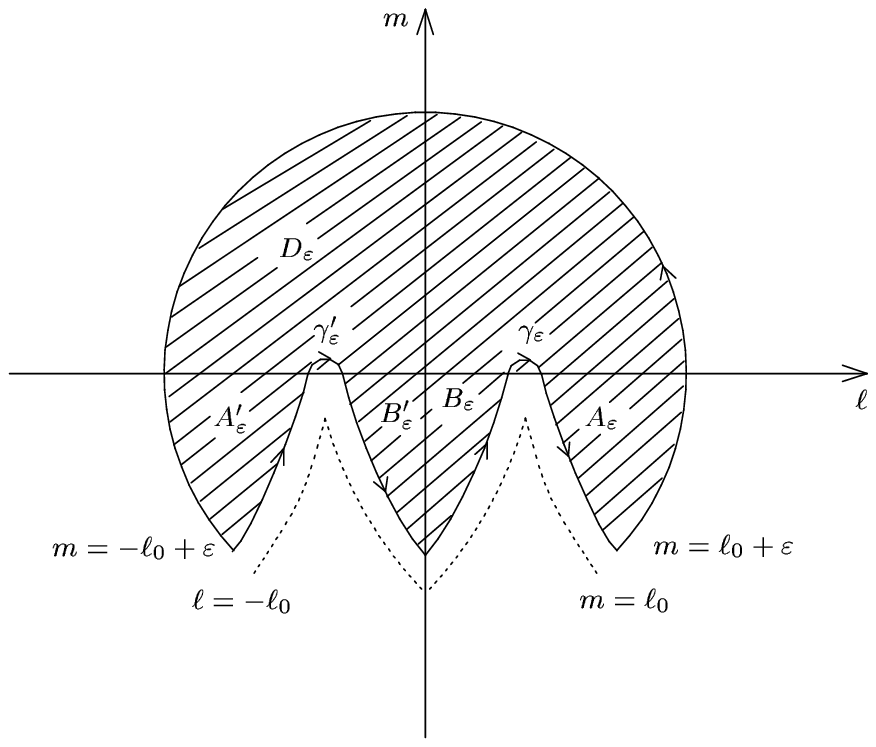, height=7.5cm, width=9cm}
\end{center}
\begin{center}
Figure 5
\end{center}
disk centered at
the origin and with radius sufficiently large so that it contains the points $(-a,0),$ $(a,0),$ and $(0,b)$
and the the support of $\vp.$ Let $D_{II,\epsilon}$ be the region obtained from $D_{II}$ by shifting it by
$\epsilon$ along the characteristics and by removing two small half disks centered at the points
$(-a,0),$ $(a,0),$ as shown in the Figure 5. $D_\epsilon$ is then the intersection of $D$ and
$D_{II,\epsilon}.$ 

By virtue of the Green's formula (\ref{gf1}), the double integral in the expression (\ref{eq31}) is to be
replaced by a contour integral along the oriented paths $A'_\epsilon,$ $\gamma'_\epsilon,$
$B'_\epsilon,$ $B_\epsilon,$ $\gamma_\epsilon,$ and $A_\epsilon.$ 

The integration along $A_\epsilon$ is similar to the one along $\bar A_\epsilon$ calculated in Section
\ref{s4} and one can see that its limit, as $\epsilon\to 0,$ is zero. The  same is true for the integral
along $A'_\epsilon.$ 

The integral along $B_\epsilon,$ similar to the integral along $\bar B_\epsilon$ (Section \ref{s4}), tends
to $\psi(\ell_0,-\ell_0)/2^{2/3},$ as $\epsilon\to 0.$  The same thing happens with the integral along
$B'_\epsilon.$ The sum of these two limits is then $2^{1/3}\psi(\ell_0,-\ell_0)$ which is equal to
the right-hand side of (\ref{eq31}). Recall that $\psi(\ell,m)=\vp((\ell+m)/6,-((\ell-m)/4)^{2/3}).$

To complete the proof one has to show that the limits, as $\epsilon\to 0,$ of the integrals
along $\gamma'_\epsilon$ and $\gamma_\epsilon$ are both zero. Since these two contours lie in the
elliptic region of the Tricomi operator, it is more convenient to use the {\em reduced elliptic} form of the Tricomi 
operator, namely,
\begin{equation}\label{eq2}
{\cal T}_e = {\p^2\over\p x^2}+{\p^2\over\p s^2}+{1\over 3s}{\p\over\p s}.
\end{equation} 
which one obtains from equation (\ref{eq1}) via the change of variables $x=x$ and $s=2y^{3/2}/3,$ and the
corresponding Green's formula in the variables $x$ and $s$. It is then a matter of verification, which we leave to the
reader, that along both  contours the integrands remain bounded and consequently both integrals along
$\gamma'_\epsilon$ and $\gamma_\epsilon$ tend to zero. The proof is then complete. \hfill $\Box$
\skp

We observe that by virtue of Proposition \ref{pr1a}, the fundamental solution $E_{II}(x,y;0,b)$ is a complex valued in 
the upper half plane  ($y>0$) plus the region in the lower half plane outside the reflected characteristics and it is real
valued in the triangle with vertices $(-a,0),$ $(a,0),$ and $(0,b).$  Since the Tricomi operator has real coefficients, the
complex conjugate $\bar E_{II}(x,y;0,b)$ as well as the real part of $E_{II}(x,y;0,b)$ are fundamental solutions for
$\cal T.$ 

Contrary to what happened in Corollary \ref{cor1} of Section \ref{s4} where we showed how the limit of 
$E_I(x,y;0,b),$ as $(0,b)\to(0,0),$ tends to the fundamental solution $F_-(x,y)$  given by formulas (\ref{eq1b}) and
(\ref{eq1b'}), neither of the fundamental solutions just obtained tend to the fundamental solution $F_+(x,y)$ given by
formulas (\ref{eq1a}) and (\ref{eq1a'}). For this to happen we need a particular linear combination of $E_{II}$ and 
$\bar E_{II}.$ Let $\lambda$ and $\mu$ be such that
$$ 
\left\{ \begin{array}{ll}
\lambda e^{i\pi\ov6} + \mu e^{-{i\pi\ov6}} &= -1/3^{1/2}\\
\lambda\ \ \ \ + \mu &= \ \ \ 1.
\end{array}
\right.
$$
Then, the following result holds:

\begin{cor}\label{cor2}
Let $E^\sharp_{II}=\lambda E_{II}+\mu {\bar E}_{II}.$ Then $E^\sharp_{II}$ is a fundamental solution
relative to $(0,b)$ that converges in the sense of distributions to $F_+(x,y),$ as $(0,b)\to(0,0).$
\end{cor}

\end{document}